\documentclass[11pt]{amsart}

\usepackage{amsmath}
\usepackage{amssymb}
\usepackage{graphicx}

\newtheorem{theorem}{Theorem}[section]
\newtheorem{corollary}[theorem]{Corollary}
\newtheorem{lemma}[theorem]{Lemma}
\newtheorem{proposition}[theorem]{Proposition}
\theoremstyle{definition}
\newtheorem{definition}[theorem]{Definition}

\numberwithin{equation}{section}

\title[Optimal shape design]{Optimal shape design for 2D heat equations in large time}

\author[E. Tr\'elat]{Emmanuel Tr\'elat}

\address[E. Tr\'elat]{Sorbonne Universit\'es, UPMC Univ Paris 06, CNRS UMR 7598, Laboratoire Jacques-Louis Lions, F-75005, Paris, France}
\email{{\tt emmanuel.trelat@upmc.fr}}

\author[C. Zhang]{Can Zhang}
\address[C. Zhang]{[1] School of Mathematics and Statistics, Wuhan University, 430072 Wuhan, China;
[2] Sorbonne Universit\'es, UPMC Univ Paris 06, CNRS UMR 7598, Laboratoire Jacques-Louis Lions, F-75005, Paris, France}
\email{\tt zhangcansx@163.com }

\author[E. Zuazua]{Enrique Zuazua}
\address[E. Zuazua]{[1] DeustoTech, Fundaci\'on Deusto, Avda Universidades, 24, 48007, Bilbao, Basque Country, Spain;
[2] Departamento de Matem\'aticas,Universidad Aut\'onoma de Madrid, 28049 Madrid, Spain;
[3] Facultad Ingenier\'ia, Universidad de Deusto, Avda. Universidades, 24, 48007 Bilbao, Basque Country, Spain;
[4] Sorbonne Universit\'es, UPMC Univ Paris 06, CNRS UMR 7598, Laboratoire Jacques-Louis Lions,  F-75005, Paris, France}

\email{\tt enrique.zuazua@uam.es }



\begin{document}

\begin{abstract}

In this paper, we investigate the asymptotic behavior of optimal designs for the shape optimization of 2D heat equations in long time horizons. The control is the shape of the domain on which heat diffuses. The class of 2D admissible shapes is the one introduced by \u{S}ver\'{a}k in \cite{S1},  of all open subsets of a given bounded open set, whose complementary sets have a uniformly bounded number of connected components. Using a $\Gamma$-convergence approach, we establish that the parabolic optimal designs converge as the length of the time horizon tends to infinity, in the complementary Hausdorff topology, to an optimal design for  the corresponding  stationary elliptic equation.

\end{abstract}

\maketitle

\begin{center}
{\it Dedicated to Professor Viorel Barbu}
\end{center}

\section{Introduction and main result}

We consider the problem of shape optimization of the heat equation in two space dimensions in the geometric framework developed by \u{S}ver\'{a}k in \cite{S1} for the optimal design problem of elliptic 
equations, where optimization is performed in the class of admissible domains characterized essentially by the fact that their complementary sets have at most a finite prescribed number of connected components.  

More precisely, the geometric setting of the admissible domains in \cite{S1} is as follows.
Given a nonempty bounded open subset $\mathcal D\subset\mathbb R^2$ (called the \emph{design region}),  we denote by $\mathcal O$ the set of all open subsets
contained in $\mathcal D$.  Let $\omega\in \mathcal O$ be  an arbitrarily fixed nonempty open subset. For each arbitrarily fixed integer $N\in\mathbb N^+$, 
we define the set of admissible designs 
\begin{equation*}
\mathcal O_\omega^N=\Big\{\Omega\in\mathcal O\;\;|\;\; \Omega\supset\omega, \;\sharp\,\Omega^c\leq N\Big\}.
\end{equation*}
Here and in the sequel, $\Omega^c=\bar{\mathcal D}\setminus\Omega$ is the complementary subset of $\Omega$ in $\mathcal D$, and
$\sharp\,\Omega^c$ denotes the number of its connected components. 

The aim of the shape optimization problem is to find the best \emph{shape or design} of the domain within that class to minimize some
cost functional depending on the domain through the solutions of some given PDE. This subject has been extensively 
studied during the last decades both for elliptic equations and for evolution problems (see, e.g., \cite{ABFJ,AMP,BB1,ChZ1,HHS,HP1,LR,P1,SZ,S1,Tartar} and references therein).
Among the methods and techniques used to solve those shape optimization problems, calculus of variations, Hadamard shape differentiation method and homogenization theory played a central role.

Let us now describe the specific problem that we address in the present paper.

\noindent \textbf{Parabolic optimal design problem}.
Let $T>0$ be arbitrary. For any $y_0\in L^2(\mathcal D)$, any $f\in L^2(\mathcal D)$ and any $z\in H^1_0(\mathcal D)$, consider the problem of minimizing the time average performance 
\begin{equation}\label{optheat}
(P^T):\inf_{\Omega\in\mathcal O_\omega^N}J^T(\Omega)\!= \!\frac{1}{T}\!\!\int_0^T\!\!\!\int_\omega \left( |y(t,x)-z(x)|^2+|\nabla y(t,x)-\nabla z(x)|^2 \right) dxdt,
\end{equation}
where $y\in C\left([0,T];L^2(\Omega)\right)\cap L^2(0,T;H_0^1(\Omega))$ satisfies the heat  equation in the domain $\Omega$
\begin{equation}\label{heat}
\left\{
\begin{split}
&\partial_t y-\triangle y=f\;\;&\text{in}\;\;\Omega\times(0,T),\\
&y=0\;\;&\text{on}\;\;\partial\Omega\times(0,T),\\
&y(\cdot,0)=y_0\;\;&\text{in}\;\;\Omega.
\end{split}\right.
\end{equation}
Here, the control variable is the \emph{shape (or design)} $\Omega$ in which the heat equation evolves, and the heat source $f$ in the equation is assumed to be \emph{independent} of time (although more general situations in which $f$ depends on $t$ but stabilizes as $t \to +\infty$ could be treated by similar methods). The target $z=z(x) \in H^1_0(\mathcal D)$ is given and, when minimizing this functional, the goal is to steer the restriction to $\omega$ of solution of the heat equation $y$ as close as possible to $z$, by an optimal choice of the shape $\Omega$ which is the domain where the Dirichlet heat equation \eqref{heat} is considered.

Since the domain $\Omega\in\mathcal O_\omega^N$ is the unknown in the above minimization problem, it is useful to note that, for the heat equation \eqref{heat} to be well posed in the functional space $C\left([0,T];L^2(\Omega)\right)\cap L^2(0,T;H_0^1(\Omega))$, it suffices that $\Omega$ be an open bounded subset of $\mathbb{R}^2$, not being necessarily of class $C^2$ (when $\Omega$ is $C^2$, we have moreover $y(t,\cdot)\in H^2(\Omega)$ for a.e. $t>0$).

We will prove further that $(P^T)$ has at least one minimizer $\Omega_T\in\mathcal O_\omega^N$.

In this  problem, the time horizon $T$ is regarded as a parameter. 
In order to investigate the long-time behavior of optimal designs for the problem $(P^T)$ as $T\to+\infty$, we next consider a reference  elliptic optimal design problem.

\noindent\textbf{Associated elliptic optimal design problem.}
For the same $z\in H^1(\omega)$ and $f\in L^2(\mathcal D)$ as above, we consider the   shape optimization problem 
\begin{equation}\label{optelliptic}
(P^s):\;\;\;\;\inf_{\Omega\in\mathcal O_\omega^N}J^s(\Omega)= \int_\omega \left( |p(x)-z(x)|^2+|\nabla p(x)-\nabla z(x)|^2 \right) dx,
\end{equation}
where $p\in H_0^1(\Omega)$ is the unique solution to the Poisson equation in $\Omega$
\begin{equation}\label{elliptic}\left\{
\begin{split}
-\triangle p&=f\;\;&\text{in}\;\;\Omega,\\
p&=0\;\;&\text{on}\;\;\partial\Omega.
\end{split}\right.
\end{equation}
Note that the control variable here is also the shape (or design), in which the equation is fulfilled.
We will prove further that $(P^s)$ has at least one minimizer $\Omega^s\in\mathcal O_\omega^N$.

\noindent\textbf{Long-time behavior.}
For each realization of the domain $\Omega$, the solution $y(t,\cdot)$ of (\ref{heat}) converges exponentially in $H^1_0(\Omega)$ as $t\to +\infty$ towards the solution of \eqref{elliptic}.

It is then natural to conjecture that the optimal shapes $\Omega_T$ for the parabolic optimal design problem (\ref{optheat}) converge (in a sense to be made precise) to optimal shapes $\Omega^s$ for the elliptic optimal design problem (\ref{optelliptic}) as $T\to +\infty$. The objective of this paper is to show that this result holds, indeed, in the geometric setting above in the complementary Hausdorff
topology  (see Section \ref{pre} for the precise definition).

In the next section, we will introduce some notations and  then briefly report on existence of minimizers for  $(P^s)$ and $(P^T)$ for  $T>0$ fixed, already established in the existing literature (see, e.g., \cite{CZ1,ChZ1,S1}).  

Numerical approximation issues for the optimal design problems  above have been addressed in \cite{CZ1,ChZ2,ChZ1}, showing that the discrete optimal shapes (defined in a finite element context)  converge in the complementary Hausdorff
topology, to an optimal shape for the continuous one as the mesh-size tends to zero.
This problem was  successfully formulated and solved in \cite{ChZ1} for   2D elliptic problems with Dirichlet boundary conditions and later extended to the heat equation case in \cite{CZ1}, and to the  wave equation  in \cite{C}.

Our objective is to address the following two specific issues:
\begin{itemize}
\item Convergence of minima: 
$$\lim_{T\rightarrow+\infty}J^T=J^s,$$
where $J^T$ and $J^s$ are  the optimal values for the problems
$(P^T)$ and $(P^s)$, respectively.
\item Convergence of minimizers:  any closure point (in complementary Hausdorff topology) as $T\rightarrow+\infty$ of minimizers of $(P^T)$ is a minimizer of $(P^s)$.
\end{itemize}

Our main result hereafter solves these two questions. 

\begin{theorem}\label{shapeconv}
Given any $y_0\in L^2(\mathcal D)$, any $f\in L^2(\mathcal D)$ and any $z\in H^1_0(\Omega)$, there exists $C>0$ (not depending on the time horizon $T$) such that
\begin{equation}\label{22519}
\big|J^T-J^s\big|\leq \frac{C}{\sqrt T}\qquad\forall T>0.
\end{equation}
Moreover, the problems $(P^s)$ and $(P^T)$, for every $T$, have at least one minimizer, and any closure point (in complementary Hausdorff topology) of minimizers of $(P^T)$ as $T\to+\infty$ is a minimizer of $(P^s)$.  
\end{theorem}

In practical applications, optimal shapes are often computed on the basis of the steady-state model, but they are then employed as quasi-optima for the time-evolving problem, often without rigorous proofs (see, e.g., \cite{Ba}). This approximation is based on the intuitive idea that, if the time-evolving dynamics converges for long time to the steady state one, elliptic optimal shapes should be nearly optimal for the time-evolution problem as well.  
From \eqref{22519} and \eqref{3081} (in the proof of Theorem~\ref{shapeconv}  below) we also derive the following result which justifies such an approximation.

\begin{corollary}
For any minimizer $\Omega^s$ of $(P^s)$, we have
$$
\big|J^T(\Omega^s)-J^T|\leq O\left(\frac{1}{\sqrt T}\right)\qquad\forall T>0.
$$
\end{corollary}

Similar results have been established in various contexts. For instance, in \cite{HS},   a shape optimization problem for the heat equation was considered, in which the support of a Radon measure on the lateral boundary was selected in an optimal way. Under certain compact assumptions, they first showed the existence of an optimal solution for this optimization problem.  They also proved convergence to an optimal solution of the corresponding stationary optimization problem for long time horizons.
Recently, the authors of \cite{AMP} have investigated the long-time behavior of a two-phase optimal design problem. More precisely,  
they considered an optimal design problem of minimizing the time average of the dissipated thermal energy during a fixed time interval $[0,T]$ and in a fixed bounded domain, where the dissipation is governed by a two-phase isotropic transient heat equation, the time \emph{independent} material properties being the design variables. 
Via a  $\Gamma$-convergence technique and the exponential decay of the energy for the heat equation,   they proved that the optimal solutions of an associated relaxed design problem converge, as $T\to+\infty$, to an optimal relaxed design of the corresponding two-phase optimization problem for the stationary heat equation. 

There is a rich literature on the limiting asymptotic behavior of optimal control problems as the time horizon goes to infinity. This problem, as previously indicated in \cite{PZ1}, is related to the so-called \emph{turnpike property}, arising mainly in economy theory (see \cite{Gr,PZ1,TZ1,Z1,Z2}).
The work \cite{PZ1} addresses the problem of long time horizon versus steady state control in the linear setting, both for finite-dimensional models, and also PDE models, namely, the heat and the wave equations, proving that, under suitable controllability assumptions and coercivity conditions in the cost functional, optimal controls and controlled trajectories (resp., adjoint states) converge exponentially to the corresponding stationary optimal controls and states (resp., adjoint states), when the time horizon tends to infinity. This result was then extended to the more general nonlinear controlled systems \cite{GTZ, TZ1,TZZ1}, in particular to a controlled system with a time-periodic cost \cite{TZZ2}. 

Note however that the problem we address in this paper is simpler in nature since the shapes under consideration are assumed to be time-independent.

The rest of this paper is organized as follows.  In Section \ref{pre}, we recall, in particular,  the definitions of the complementary Hausdorff
topology
and the main results in $\Gamma$-convergence. Section \ref{pro} is devoted to the 
proof of  Theorem~\ref{shapeconv}.  Finally, in Section \ref{con} we conclude this paper with some further comments and open problems.

\medskip

\section{Preliminaries}\label{pre}
\subsection{Existence of optimal designs}
We first recall the definition of the Hausdorff topology and of the complementary Hausdorff topology. 
\begin{definition}
The Hausdorff  distance between 
two compacts sets $K_1$ and $K_2$ in $\mathbb R^2$ is defined by
$$
d_{H}(K_1,K_2)=\max \left( \max_{x\in K_2}\min_{y\in K_1}\|x-y\|,
\max_{x\in K_1}\min_{y\in K_2}\|x-y\|\right),
$$
where $\|\cdot\|$ is the Euclidean norm  in  $\mathbb R^2$.
\end{definition}
Recall that $\mathcal O$ is the set of all open subsets of $\mathcal D$.  
For any $\Omega_i\in \mathcal O$, $i=1,2$, we define the complementary Hausdorff distance by
\begin{equation*}
d_{H^c}(\Omega_1,\Omega_2)=\max \left(\max_{x\in\Omega^c_2}\min_{y\in\Omega_1^c}\|x-y\|,
\max_{x\in\Omega^c_1}\min_{y\in\Omega_2^c}\|x-y\|\right),
\end{equation*}
where $\Omega_i^c=\bar {\mathcal D}\setminus \Omega_i$, $i=1,2$. Then, $(\mathcal O, d_{H^c}(\cdot,\cdot))$ is a 
complete metric space.

We say that $\Omega_n\stackrel{H^c}{\longrightarrow}\Omega$ if and only if $d_{H^c}(\Omega_n,\Omega)\longrightarrow 0$, as $n\rightarrow +\infty$.

We refer  the interested reader to \cite{ChZ1} for properties related to the Hausdorff convergence and facts that might seem counterintuitive. For example, the convergence of $\{\Omega_n\}_{n\geq1}$ to $\Omega$ in the $H^c$ topology does not guarantee the convergence of the Lebesgue measure of $\Omega_n$ to that of  $\Omega$.

For each fixed $N$ and each open subset $\omega$, the set of admissible designs $\mathcal O_\omega^N$,  as defined in the introduction, 
 is well known  (see, e.g., \cite{CZ1,ChZ1,S1})  to be compact for the complementary Hausdorff topology. This implies that, for any sequence $(\Omega_j)_{j\geq1}$ of $ \mathcal O_\omega^N$, there exist $\Omega \in \mathcal O_\omega^N$ and a subsequence $\{\Omega_k\}_{k\geq1}$ of $\{\Omega_j\}_{j\geq1}$ such that $\Omega_k\stackrel{H^c}{\longrightarrow} \Omega$ as $k\rightarrow+\infty$.  

For any $\Omega\subset\mathcal O$, $H_0^1(\Omega)$ is  defined as the closure, for the $H_0^1(\Omega)$ topology, of all smooth functions with compact support in $\Omega$. Accordingly, any function of $H_0^1(\Omega)$ can be extended by 0 to a function of $H_0^1(\mathbb R^2)$ (and $H_0^1(\mathcal D)$).
Here and in the sequel, for any $y\in H_0^1(\Omega)$ with  $\Omega\in\mathcal O$,  we will denote by $\widetilde y$ 
its extension by zero to the fixed domain $\mathcal D$. 

Next, we introduce the notion of
 $\Gamma$-convergence for open subsets, which plays a crucial role in the investigation of existence of optimal designs in shape optimization problems. 

\begin{definition}
We say that $\Omega_n\stackrel{\Gamma}{\longrightarrow}\Omega$ if for any $f\in L^2(\mathcal D)$, the solution $p_n$  of the Poisson equation
\begin{equation*}\left\{
\begin{split}
-\triangle p_n&=f\;\;&\text{in}\;\;\Omega_n,\\
p_n&=0\;\;&\text{on}\;\;\partial\Omega_n,
\end{split}\right.
\end{equation*}
satisfies
$$\widetilde{p}_n\longrightarrow \widetilde p\;\;\;\text{ in} \;\;\; H_0^1(\mathcal D),$$ 
where $p$ is the solution to
\begin{equation*}\left\{
\begin{split}
-\triangle p&=f\;\;&\text{in}\;\;\Omega,\\
p&=0\;\;&\text{on}\;\;\partial\Omega.
\end{split}\right.
\end{equation*}
\end{definition}

In general, $H^c$-convergence does not imply $\Gamma$-convergence. Indeed, it is well known that homogenization phenomena may occur at the limit, when the sequence of designs is allowed to develop an increasing number of holes. In this case the limit of the solutions of the Dirichlet-Laplacian may be the solution of a different elliptic problem (see, e.g., \cite{ABFJ,Tartar}). Fortunately, several situations are known where the $H^c$-convergence does imply the $\Gamma$-convergence and the above relaxation  phenomena do not occur (see, e.g., \cite[Theorem 4.6.7]{BB1}). The following one is due to V. \u{S}ver\'{a}k.

\begin{theorem}[\cite{S1}]\label{prop2}
Let $\Omega$ and $(\Omega_n)_{n\geq1}$ belong to $\mathcal O_\omega^N$. Then $\Omega_n\stackrel{\Gamma}{\longrightarrow}\Omega$ is equivalent to $\Omega_n\stackrel{H^c}{\longrightarrow}\Omega$.
\end{theorem}

Since $\mathcal O_\omega^N$ is  compact in the complementary Hausdorff topology, from
Theorem~\ref{prop2} we deduce  that for any sequence
of designs $(\Omega_n)_{n\geq1}\subset\mathcal O_\omega^N$, there exist $\Omega\in\mathcal O_\omega^N$ and a subsequence (for simplicity we still denote it in the same way), such that $\Omega_n\stackrel{H^c}{\longrightarrow}\Omega$ and
$\Omega_n\stackrel{\Gamma}{\longrightarrow}\Omega$.

As corollaries of Theorem~\ref{prop2}, the existence of minimizers of the shape optimization problem $(P^s)$, as well as $(P^T)$ with each $T>0$, have already been established. We now  state it as follows.
\begin{proposition}\label{bao3}
The problem $(P^s)$ has at least one minimizer, and for any $T>0$, the problem $(P^T)$ has at least one minimizer.
\end{proposition}

For a proof, we refer the interested reader to \cite{S1} or \cite{ChZ1} for the elliptic optimal design problem, and to \cite{CZ1} for the heat one. 

Uniqueness of optimal solutions is still an open and challenging  issue in the theory of shape optimization problems. For example, the authors of \cite{ABFJ} constructed a specific example for which there is an infinite number of optimal designs.

\subsection{The uniform Poincar\'e inequality} 
We recall that for each open subset $\Omega\in\mathcal O$, the first eigenvalue $\lambda_1(\Omega)$ for the Laplace operator $-\triangle$ in $\Omega$, with zero Dirichlet boundary conditions, is given by the Rayleigh formula (see, e.g., \cite{Brezis})
\begin{equation*}
\lambda_1(\Omega)=\inf_{u\in H_0^1(\Omega)\setminus\{0\}}\frac{\int_\Omega |\nabla u(x)|^2\,dx}{\int_\Omega|u(x)|^2\,dx}.
\end{equation*}
Minimization problems for elliptic eigenvalue problems have received significant attention in the literature since
the first result by Faber and Krahn, concerning the  first eigenvalue of the Laplace operator $-\triangle$ in 2D, with Dirichlet boundary conditions, among open subsets with equal area, ensuring that $\lambda_1(\Omega)\geq\lambda_1(B)>0$ for every $\Omega\in\mathcal O$, where $B$ is a ball in $\mathbb R^2$ with area equal to the Lebesgue measure of $\mathcal D$ (see, e.g., \cite[Chapter 6]{BB1} or \cite{HP1}).
Consequently, the following Poincar\'e inequality holds uniformly in the class of open sets $\mathcal O$, which will play a crucial role in the proof of Theorem~\ref{shapeconv}.

\begin{lemma}\label{lemma3}
There exists $C>0$ depending only on the area of $\mathcal D$,  such that 
\begin{equation}\label{unip}
\int_\Omega |u(x)|^2\,dx\leq C \int_\Omega |\nabla u(x)|^2\,dx,
\end{equation}
 for all $\Omega\in \mathcal O$ and $u\in H_0^1(\Omega)$. 
\end{lemma}

\section{Proof of Theorem~\ref{shapeconv}}\label{pro}

From Proposition~\ref{bao3}, we have seen that the shape optimization problems $(P^s)$ and $(P^T)$, for any fixed $T>0$, have minimizers in the class of admissible shapes $\mathcal O_\omega^N$. Based on a  $\Gamma$-convergence argument, we next prove the long-time behavior of the optimal design problems $(P^T)$ stated in Theorem~\ref{shapeconv}.  For an introduction to the theory of $\Gamma$-convergence in the calculus of variations, the interested reader is referred to \cite{D}.

\begin{proof}[\textbf{Proof of Theorem~\ref{shapeconv}}.]
We proceed in three steps.

\noindent\textbf{Step 1.} We first show the upper bound 
\begin{equation}\label{2261}
J^T-J^s\leq \frac{C}{\sqrt T}\qquad\forall T>0,
\end{equation}
for some constant $C=C(|\mathcal D|, \|f\|_{L^2(\mathcal D)},\|y_0\|_{L^2(\mathcal D)},\| z\|_{H^1_0(\mathcal D)})>0$ not depending on $T$. Recall that $J^T$ and $J^s$ are, respectively,  the optimal values for the problems
$(P^T)$ and $(P^s)$.

Assume that $\Omega^s\in\mathcal O^N_\omega$ is an optimal design of $(P^s)$. Then $J^s=J^s(\Omega^s)$. Since $\Omega^s$ is an admissible design of $(P^T)$,  we obviously have $J^T\leq J^T(\Omega^s)$.
Hence
\begin{equation}\label{2262}
J^T-J^s\leq J^T(\Omega^s)-J^s(\Omega^s).
\end{equation}
Now, let us assume that $p^s\in H_0^1(\Omega^s)$ is the solution of
\begin{equation}\label{jingbao1}
\left\{
\begin{split}
-\triangle p^s&=f\;\;&\text{in}\;\;\Omega^s,\\
p^s&=0\;\;&\text{on}\;\;\partial\Omega^s.
\end{split}\right.
\end{equation}
The energy identity ensures that
$\int_{\Omega^s} |\nabla p^s(x)|^2\,dx=\int_{\Omega^s}f(x)p^s(x)\,dx$.
By the uniform Poincar\'e inequality \eqref{unip} in Lemma~\ref{lemma3} and the Cauchy-Schwarz inequality, there exists $C=C(|\mathcal D|)>0$ such that
\begin{equation}\label{2246}
\|\nabla p^s\|_{L^2(\Omega^s)}\leq C\|f\|_{L^2(\mathcal D)} ,
\end{equation}
and
\begin{equation}\label{baojian1}
\|p^s\|_{L^2(\Omega^s)}\leq C\|f\|_{L^2(\mathcal D)}.
\end{equation}
We denote by $y_T(\cdot)\in L^\infty(0,T;L^2(\Omega^s))\cap L^2(0,T;H_0^1(\Omega^s))$ the solution of
\begin{equation}\label{2241}\left\{
\begin{split}
&\partial_t y_T-\triangle y_T=f(x)\;\;&\text{in}\;\;\Omega^s\times(0,T),\\
&y_T=0\;\;\;&\text{on}\;\;\partial \Omega^s\times(0,T),\\
&y_T(x,0)=y_0\;\;&\text{in}\;\;\Omega^s.
\end{split}\right.
\end{equation}
Let $(S(t))_{t\geq0}$ be the $C_0$ semigroup in $L^2(\Omega^s)$ generated by the Laplace operator $\triangle$ on the domain $D(\triangle)=\{u\in H_0^1(\Omega^s)\ \mid\ \triangle u \in L^2(\Omega^s)\}$ (see, e.g., \cite{Marius}).
Energy estimates ensure that $\|S(t)\|_{\mathcal L(L^2(\Omega^s);L^2(\Omega^s))}\leq e^{-\lambda t}$ for every $t\geq0$, with $\lambda=\lambda(|\mathcal D|)>0$.
Since $y_T(t)=S(t)y_0+\int_0^tS(t-\tau)f\,d\tau$ for every $t\in[0,T]$, we infer that
\begin{equation}\label{22421}
\max_{t\in[0,T]} \|y_T(t)\|_{L^2(\Omega^s)}\leq \|y_0\|_{L^2(\mathcal D)}+\frac{1}{\lambda}\|f\|_{L^2(\mathcal D)}.
\end{equation}
Multiplying by $y_T(\cdot)$ the equation \eqref{2241} and integrating by parts, we get
\begin{multline*}
\frac{1}{2}\left(\|y_T(T)\|^2_{L^2(\Omega^s)}-\|y_0\|^2_{L^2(\Omega^s)}\right)+\int_0^T\|\nabla y_T(t)\|_{L^2(\Omega^s)}^2\,dt\\
\leq \int_0^T\int_{\Omega^s}f(x)y_T(x,t)\,dx\, dt
\leq T\|f\|_{L^2(\mathcal D)}\max_{t\in[0,T]} \|y_T(t)\|_{L^2(\Omega^s)}.
\end{multline*}
Combined with \eqref{22421}, this implies that
\begin{equation}\label{22410}
\int_0^T\|\nabla y_T(t)\|_{L^2(\Omega^s)}^2\,dt \leq CT ,
\end{equation}
for some constant $C=C(|\mathcal D|, \|f\|_{L^2(\mathcal D)},\|y_0\|_{L^2(\mathcal D)})>0$ not depending on $T$.

Next, we set $\delta y_T(t)= y_T(t)-p^s$, for every $t\in[0,T]$. 
It follows from \eqref{jingbao1} and \eqref{2241} that $\delta y_T(\cdot)$ is solution of the heat equation in $\Omega^s$,
\begin{equation*}\left\{
\begin{split}
&\partial_t \delta y_T-\triangle \delta y_T=0\;\;&\text{in}\;\;\Omega^s\times(0,T),\\
&\delta y_T=0\;\;\;&\text{on}\;\;\partial \Omega^s\times(0,T),\\
&\delta y_T(x,0)=y_0-p^s\;\;&\text{in}\;\;\Omega^s.
\end{split}\right.
\end{equation*}
It is easy to see that there exists $C(|\mathcal D|)>0$ (not depending on $T$) such that
$\int_0^T\|\delta y_T(t)\|^2_{L^2(\Omega^s)}\,dt\leq C(|\mathcal D|)\|y_0-p^s\|^2_{L^2(\Omega^s)}$
and $\int_0^T\|\nabla \delta y_T(t)\|_{L^2(\Omega^s)}^2\,dt \leq \|y_0-p^s\|^2_{L^2(\Omega^s)}$.
These last two inequalities, combined with  \eqref{2246}, imply that
\begin{equation}\label{2249}
\int_0^T\|\delta y_T(t)\|^2_{L^2(\Omega^s)}\,dt+\int_0^T\|\nabla \delta y_T(t)\|_{L^2(\Omega^s)}^2\,dt \leq C ,
\end{equation}
for some constant $C=C(|\mathcal D|, \|f\|_{L^2(\mathcal D)},\|y_0\|_{L^2(\mathcal D)})>0$
not depending on $T$.

Note that 
\begin{equation}\label{2265}
J^T(\Omega^s)-J^s(\Omega^s)=I_1+I_2
\end{equation}
with
\begin{multline*}
I_1= \frac{1}{T}\int_0^T\left(\|y_T(t)-z\|_{L^2(\omega)}+\|p^s-z\|_{L^2(\omega)}\right)\\
\times
\left(\|y_T(t)-z\|_{L^2(\omega)}-\|p^s-z\|_{L^2(\omega)}\right)\,dt
\end{multline*}
and
\begin{multline*}
I_2=\frac{1}{T}
\int_0^T\left(\|\nabla y_T(t)-\nabla z\|_{L^2(\omega)}+\|\nabla p^s-\nabla z\|_{L^2(\omega)}\right)\\
\times
\left(\|\nabla y_T(t)-\nabla z\|_{L^2(\omega)}-\|\nabla p^s-\nabla z\|_{L^2(\omega)}\right)\,dt.
\end{multline*}
We first estimate the term $I_1$ as follows. By the triangle inequality
$$
\Big|\|y_T(t)-z\|_{L^2(\omega)}-\|p^s-z\|_{L^2(\omega)}\Big|\leq \|y_T(t)-p^s\|_{L^2(\omega)}\leq \|\delta y_T(t)\|_{L^2(\Omega^s)},
$$
 we get that
\begin{multline*}
|I_1|\leq 
\frac{1}{T}\left(\int_0^T\left(\max_{t\in[0,T]}\|y_T(t)\|_{L^2(\Omega^s)}+\|p^s\|_{L^2(\Omega^s)}+2\|z\|_{L^2(\omega)}
\right)^2dt\right)^{1/2}\\
\times
\left(\int_0^T\|\delta y_T(t)\|^2_{L^2(\Omega^s)}\,dt\right)^{1/2}.
\end{multline*}
This, together with \eqref{baojian1}, \eqref{22421} and \eqref{2249}, leads to 
\begin{equation}\label{2264}
|I_1|\leq \frac{C}{\sqrt T} ,
\end{equation}
for some constant $C=C(|\mathcal D|, \|f\|_{L^2(\mathcal D)},\|y_0\|_{L^2(\mathcal D)},\|z\|_{L^2(\omega)})>0$.
Similarly, the term $I_2$ is estimated by
\begin{multline*}
|I_2|\leq \frac{1}{T}\left(\int_0^T\big(\|\nabla y_T(t)\|_{L^2(\omega)}+\|\nabla p^s\|_{L^2(\omega)}+2\|\nabla z\|_{L^2(\omega)}\big)^2\,dt\right)^{1/2}\\   
\times
\left(\int_0^T\|\nabla \delta y_T(t)\|_{L^2(\Omega^s)}^2\,dt\right)^{1/2}.
\end{multline*}
Combined with \eqref{2246}, \eqref{22410} and \eqref{2249}, this implies that
\begin{equation}\label{2258}
|I_2|\leq \frac{C}{\sqrt T} ,
\end{equation}
for some constant $C=C(|\mathcal D|, \|f\|_{L^2(\mathcal D)},\|y_0\|_{L^2(\mathcal D)},\|\nabla z\|_{L^2(\omega)})>0$.
Therefore,  we obtain from \eqref{2265},  \eqref{2264} and \eqref{2258} that
\begin{equation}\label{3081}
\big|J^T(\Omega^s)-J^s(\Omega^s)\big|\leq\frac{C}{\sqrt T},
\end{equation}
with $C>0$ as above (not depending on $T$).
The estimate \eqref{2261} now follows from \eqref{2262} and \eqref{3081}.

\noindent\textbf{Step 2.}
Let us establish the lower estimate
\begin{equation}\label{225goal}
J^T-J^s\geq - \frac{C}{\sqrt T}\qquad\forall T>0 ,
\end{equation}
for some constant $C=C(|\mathcal D|, \|y_0\|_{L^2(\mathcal D)}, \|f\|_{L^2(\mathcal D)},\|z\|_{H^1(\omega)})>0$ not depending on $T$.

For any $T>0$, we assume that $\Omega^T\in\mathcal O_\omega^N$ is a minimizer of $(P^T)$. Reasoning as in \eqref{2262}, we have
\begin{equation}\label{22611}
J^T-J^s\geq J^T(\Omega^T)-J^s(\Omega^T).
\end{equation}
Let $y^T\in L^\infty\big(0,T;L^2(\Omega^T)\big)\cap L^2\big(0,T;H_0^1(\Omega^T)\big)$ be the corresponding solution to the optimal design
$\Omega^T$ for the problem $(P^T)$. Using the arguments employed to obtain the estimates \eqref{22421} and \eqref{22410}, we also have that
\begin{equation}\label{2253}
\max_{t\in[0,T]}\|y^T(t)\|_{L^2(\Omega^T)}\leq C ,
\end{equation} 
and 
\begin{equation*}
\int_0^T\|\nabla y^T(t)\|_{L^2(\Omega^T)}^2\,dt\leq CT ,
\end{equation*}
for some constant $C=C(|\mathcal D|, \|f\|_{L^2(\mathcal D)},\|y_0\|_{L^2(\mathcal D)})>0$ not depending on $T$.

Now, let $p_T\in H_0^1(\Omega^T)$ be the solution of
\begin{equation*}\left\{
\begin{split}
-\triangle p_T=f\;\;\;&\text{in}\;\;\Omega^T,\\
p_T=0\;\;\;&\text{on}\;\;\partial\Omega^T.
\end{split}\right.
\end{equation*}
By the uniform Poincar\'e inequality \eqref{unip}, 
as in \eqref{2246} and \eqref{baojian1},
there exists $C=C(|\mathcal D |)>0$ (not depending on $T$) such that
\begin{equation*}
\|\nabla p_T\|_{L^2(\Omega^T)}\leq C\|f\|_{L^2(\mathcal D)} ,
\end{equation*}
and 
\begin{equation}\label{2254}
\|p_T\|_{L^2(\Omega^T)}\leq C\|f\|_{L^2(\mathcal D)}.
\end{equation}

Then $\delta y^T(t)= y^T(t)-p_T$, $t\in[0,T]$ is solution of the heat equation on $\Omega^T$. 
Reasoning as in \eqref{2249}, we obtain that 
\begin{equation}\label{2255}
\int_0^T\|\delta y^T(t)\|^2_{L^2(\Omega^T)}\,dt+\int_0^T\|\nabla \delta y^T(t)\|_{L^2(\Omega^T)}^2\,dt
\leq C \qquad\forall T>0,
\end{equation}
for some constant $C=C(|\mathcal D|, \|y_0\|_{L^2(\mathcal D)}, \|f\|_{L^2(\mathcal D)})>0$.

We now write $J^T(\Omega^T)-J^s(\Omega^T)=I_3+I_4$ with
\begin{multline}\label{2251}
I_3= \frac{1}{T}\int_0^T\left(\|y^T(t)-z\|_{L^2(\omega)}+\|p_T-z\|_{L^2(\omega)}\right)\\
\times
\left(\|y^T(t)-z\|_{L^2(\omega)}-\|p_T-z\|_{L^2(\omega)}\right)\,dt ,
\end{multline}
and
\begin{multline*}
I_4= \frac{1}{T}
\int_0^T\left(\|\nabla y^T(t)-\nabla z\|_{L^2(\omega)}+\|\nabla p_T-\nabla z\|_{L^2(\omega)}\right)\\
\times
\left(\|\nabla y^T(t)-\nabla z\|_{L^2(\omega)}-\|\nabla p_T-\nabla z\|_{L^2(\omega)}\right)\,dt.
\end{multline*}
For any $t\in(0,T)$, by the triangle inequality,
\begin{equation*}
    \Big|\|y^T(t)-z\|_{L^2(\omega)}-\|p_T-z\|_{L^2(\omega)}\Big|\leq
    \|y^T(t)-p_T\|_{L^2(\omega)}\leq 
     \|\delta y^T(t)\|_{L^2(\Omega^T)},
\end{equation*}
and from \eqref{2251} we see  that
\begin{multline}\label{22510}
|I_3|\leq
\frac{1}{T}\left(\int_0^T\left(\max_{t\in[0,T]}\|y^T(t)\|_{L^2(\Omega^T)}+\|p_T\|_{L^2(\Omega^T)}+2\|z\|_{L^2(\omega)}\right)^2dt\right)^{1/2}\\
\times
\left(\int_0^T\|\delta y^T(t)\|^2_{L^2(\Omega^T)}\,dt\right)^{1/2}.
\end{multline}
This, together with \eqref{2253}, \eqref{2254} and \eqref{2255}, implies that $|I_3|\leq C/\sqrt T$ for some constant $C=C(|\mathcal D|, \|y_0\|_{L^2(\mathcal D)}, \|f\|_{L^2(\mathcal D)},\|z\|_{L^2(\omega)})$ not depending on $T$. Also, similar arguments as those for \eqref{2258}
lead to 
\begin{equation}\label{2259}
|I_4|\leq \frac{C}{\sqrt T}
\end{equation}
for some constant $C=C(|\mathcal D|, \|y_0\|_{L^2(\mathcal D)}, \|f\|_{L^2(\mathcal D)},\|\nabla z\|_{L^2(\omega)})>0$.

Hence, it follows from \eqref{22510}  and \eqref{2259} that
\begin{equation}\label{22534}
\big|J^T(\Omega^T)-J^s(\Omega^T)\big|\leq \frac{C}{\sqrt T}\qquad\forall T>0 ,
\end{equation}
for some constant $C=C(|\mathcal D|, \|y_0\|_{L^2(\mathcal D)}, \|f\|_{L^2(\mathcal D)},\| z\|_{H^1(\omega)})>0$. Combined with \eqref{22611}, this implies \eqref{225goal}. 

From Steps 1 and 2, the estimate \eqref{22519} is now established. 

\noindent\textbf{Step 3.}
Finally, we now establish the long-time behavior of minimizers of $(P^T)$.  Let $(T_n)_{n\geq1}$ be an increasing sequence of positive times such that $\lim_{n\rightarrow+\infty}T_n=+\infty$. For each $T_n$, we assume that $\Omega^{T_n}\in\mathcal O^N_\omega$ is an optimal design for  $(P^{T_n})$. Since $\mathcal O^N_\omega$ is compact in the complementary Hausdorff topology, up to a subsequence (still denoted with the same notation), there exists $\Omega^*\in\mathcal O^N_\omega$ such that
$\Omega^{T_n}\stackrel{H^c}{\longrightarrow} \Omega^*$. Our goal  is to show that $\Omega^*$ is an optimal design for  $(P^s)$, i.e., that $J^s( \Omega^*)=J^s$.

 From Theorem~\ref{prop2}, we have $\Omega^{T_n}\stackrel{\gamma}{\longrightarrow} \Omega^*$. In other words, the solution $p_n\in H_0^1(\Omega^{T_n})$ of
\begin{equation*}  
\left\{
\begin{split}
-\triangle p_n&=f\;\;&\text{in}\;\;\Omega^{T_n},\\
p_n&=0\;\;&\text{on}\;\;\partial \Omega^{T_n}
\end{split}\right.
\end{equation*}
satisfies $\widetilde p_n\longrightarrow\widetilde p^*\;\;\;\text{in}\;\;H_0^1(\mathcal D)\;\;\text{as}\;\;n\rightarrow+ \infty$,
where $p^*\in H_0^1(\Omega^*)$ is the solution of
\begin{equation*}  
\left\{
\begin{split}
-\triangle p^*&=f\;\;&\text{in}\;\;\Omega^{*},\\
p^*&=0\;\;&\text{on}\;\;\partial \Omega^{*}
\end{split}\right.
\end{equation*}
Hence,
\begin{equation}\label{22535}
J^s(\Omega^{T_n})\rightarrow J^s(\Omega^*)\;\;\text{as}\;\;n\rightarrow +\infty.
\end{equation}
Note that, for any $n\in\mathbb N$, we have $J^{T_n}=J^{T_n}(\Omega^{T_n})$ and 
$$
\big|J^s-J^s(\Omega^*)\big|\leq\big| J^s-J^{T_n}\big|+\big| J^{T_n}(\Omega^{T_n})-J^s(\Omega^{T_n})\big|+\big|J^s(\Omega^{T_n})-J^s(\Omega^*)\big|.
$$
By letting $n$ tend to infinity in the above inequality, we get from \eqref{22519}, \eqref{22534} and \eqref{22535} that $J^s(\Omega^*)=J^s$. This completes the proof.
\end{proof}

\section{Conclusions and further comments}\label{con}
In this paper, we have established by $\Gamma$-convergence techniques that the optimal designs for heat equations converge, as the time horizon tends to infinity, towards an optimal design of the corresponding design problem for the elliptic Poisson equation, in the sense of complementary Hausdorff topology.  

Several remarks are in order.

\noindent\textbf{More general operators.}
In this paper, we fully rely on the geometric setting of admissible designs and on the results established in \cite{S1}, and therefore our convergence result is restricted to 2D Dirichlet problem.  Although we only considered the Dirichlet-Laplacian operator, by the same techniques, it is likely that the results of this paper also hold for more general 2D elliptic operators in divergence form with Dirichlet boundary conditions, and for 2D elliptic Stokes system with Dirichlet boundary conditions. We refer the reader to  \cite{BB1,ChZ1,S1}, for instance, for a discussion of elliptic optimal design problems for those models. 

\noindent\textbf{Higher dimension.}
The method developed here may certainly be adapted to deal with the heat equation in higher dimension, in an appropriate class of admissible domains. Note that the proof of our main result relies on the following two key facts: 
\begin{itemize}
\item[(i)] Compactness of admissible domains in the complementary Hausdorff topology.
Compactness holds in higher dimension in more restricted classes of domains obtained,  for instance, by imposing uniform $BV$-norm of the boundaries, on the perimeter, or by imposing the uniform exterior cone property (see, e.g.,  \cite{BB1} and \cite[Page 1083]{PTZ1}).
\item[(ii)] The $\Gamma$-convergence property of domains, allowing to pass to the limit on the solutions of the Dirichlet elliptic problem. It can be guaranteed to hold, for instance,  in the class of convex sets, the class of domains satisfying a uniform exterior cone property, or the class of domains satisfying a uniform capacity density condition (see \cite[Theorem 4.6.7]{BB1}).
\end{itemize}

\noindent\textbf{Damped wave equation.}
Our results and proofs heavily rely on the exponential decay of the energy for the heat equation in a given domain. Accordingly, our methods also apply for the shape optimization of strongly damped wave equations in the geometric setting by \u{S}ver\'{a}k (see \cite{C} for the extension of results in \cite{S1} to the wave equation). 

However, because of the lack of exponential decay for the conservative Schr\"odinger and wave equations, the long-time behavior of shape optimization for these two equations is an open problem. In fact, for conservative problems, it could well be that the optimal shapes $\Omega_T$ reproduce the oscillatory pattern of solutions as $T$ increases.

\noindent\textbf{Time-dependent source term.}
The right-hand side term $f$ has been taken to be independent of $t$. But, as mentioned in the introduction, one could also consider time-varying forcing terms $f=f(t,x)$ under the condition that they converge exponentially to a steady applied force $f^*$ as $t\to+\infty$.

\noindent\textbf{Convergence rates.}
We proved that the optimal designs  for heat equations converge, as the time horizon tends to infinity, towards an optimal one for the stationary heat equation.  
Obtaining convergence rates is of interest, but this subject is completely open. 

This issue is even open for simpler problems. 
For instance,  in \cite{PZ2}, an optimal control problem in a fixed domain, with an applied right-hand side time-independent forcing control, was considered for a semilinear heat equation.
By $\Gamma$-convergence arguments, optimal controls were proved to converge to the steady-state ones. But convergence rates have not been derived. 

Using the optimality system and Linear Quadratic Riccati theory, by means of perturbations arguments,  convergence rates were proved under suitable smallness conditions on the target for semilinear heat equations.  
Optimality conditions could also be useful in the context of shape optimization. 
But they usually require a more limited geometric setting so that Hadamard shape derivatives can be employed (see, e.g., \cite{HP1,S11,asz,SZ}). Whether this suffices to achieve convergence rates is an interesting open problem.

\noindent\textbf{Shape turnpike.}
In the context of time-varying shapes, the turnpike problem is completely open (see \cite{TZ1,TZZ1,TZZ2}). The possible stabilization of optimal designs in large time, when allowing the design to evolve in time as well as the evolution problem, is a much more complex problem than the one we addressed here.

\noindent\textbf{Initial data fixed or not.}
We have worked with fixed initial data and right-hand side terms but one could consider more general situations. For instance, there are at least two possible ways to allow the initial data to vary:
\begin{itemize}
\item[(i)] Initial data depend (only) on the time horizon $T$ and are all bounded uniformly. 
\item[(ii)] Robust optimal shape designs: Initial data vary, for instance, in the unit ball $B_1$ of $L^2(\mathcal D)$. One can then define the optimal design problem in some uniform manner with respect to all these initial data, by considering the min-max cost
$$
\min_{\Omega\in\mathcal O_\omega^N}\max_{y_0\in B_1}J^T(\Omega,y_0).
$$
\end{itemize}
Since the constants in the proof of Theorem \ref{shapeconv} depend on the $L^2$-norm of the initial data, the method of this paper can be applied to handle these problems.

It would also be interesting to consider shape optimization problems (and turnpike issues) for the heat equation with random initial data.

Randomization has been shown to be a useful tool for a number of optimal shape design problems (see \cite{PTZ2,PTZ1,PTZ3}). In these works the PDE was formulated on a fixed reference domain and the shapes to be optimized were the location of sensors and actuators.
Through randomization, the average value of the cost functional turned out to have a spectrally diagonal structure.    
But, in these papers, the fact that the PDE under consideration was settled on a fixed domain $\Omega$ played an important role, since this allowed the randomization procedure to be defined in the basis of eigenfunctions of the Dirichlet-Laplacian on this domain. However, in the present context, the domain where the PDE holds varies, being the control variable. The way of randomization needs to be implemented so as to simplify the cost under consideration is an interesting open problem.

\noindent\textbf{Terminal constraints.}
We have treated  the shape optimization problem by letting the terminal state $y(T)$ free.
It would be interesting to address similar problems in the context of controllability, the goal being to drive the solution to some given target, employing time-varying shapes $t \rightarrow \Omega(t)$ as controls. The problem of controlling the wave and Schr\"odinger equations using the shape of the domain as control parameter has been analyzed, for instance, in \cite{Beau, Gugat1, Ivan}.
There is plenty of issues to investigate in that setting for heat-like equations and in particular to investigate the turnpike property.

\bigskip

\noindent \textbf{Acknowledgment}.
The authors acknowledge the financial support by the grant FA9550-14-1-0214 of the EOARD-AFOSR. The second author was partially supported by the National Natural Science Foundation of China under grants 11501424 and 11371285.
The third author was partially supported by the Advanced Grant DYCON (Dynamic Control) of the European Research Council Executive Agency, FA9550-15-1-0027 of AFOSR, the MTM2014-52347 Grant of the MINECO (Spain) and ICON of the French ANR.


\begin{thebibliography}{99}

\bibitem{ABFJ}
G. Allaire, E. Bonnetier, G. Francfort,  F. Jouve, Shape optimization by the homogenization method. Numer. Math. 76 (1997),  27-68.

\bibitem{AMP}
G. Allaire, A. M\"unch, F. Periago, 
Long time behavior of a two-phase optimal design for the heat equation. 
SIAM J. Control Optim. 48 (2010),  5333-5356. 

\bibitem{Ba}
M. J. Balas,  Optimal quasi-static shape control for large aerospace antennae.  Journal of Optimization Theory and Applications,
 46 (1985), 153-170.
 
 
 
 \bibitem{Beau}
K. Beauchard, Controllability of a quantum particule in a 1D variable domain. ESAIM-COCV
14 (2008), 105-147.

\bibitem{Brezis}
H. Brezis, Functional Analysis, Sobolev Spaces and Partial Differential Equations. Springer Science \& Business Media, 2010.

\bibitem{BB1} D. Bucur, G. Buttazzo, Variational Methods in Shape Optimization Problems. Progress in Nonlinear Differential Equations 65, Birkha\"user Verlag, Basel, 2005.




\bibitem{C}
M. Cea, Optimal design for 2D wave equations. Optimization Methods \& Software, 32 (2017), 86-108.

\bibitem{CZ1} M. Cea,  E. Zuazua, 
Finite element approximation of 2D parabolic optimal design problems.  Numerical mathematics and advanced applications, 151-176, Springer, Berlin, 2006. 

\bibitem{ChZ2}
D. Chenais, E. Zuazua, 
Controllability of an elliptic equation and its finite difference. Numer. Math., 95 (2003), 63-99.

\bibitem{ChZ1} D. Chenais, E. Zuazua, Finite-element approximation of 2D elliptic optimal design. J. Math. Pures Appl.  85 (2006),  225-249. 

\bibitem{D}
G. Dal Maso,  An Introduction to $\Gamma$-Convergence. Springer Science \& Business Media, 2012.

\bibitem{DZ}
R. Dziri, J. P. Zol\'esio,  Dynamical shape control in non-cylindrical Navier-Stokes equations. Journal of Convex Analysis 6 (1999),293-318. 

\bibitem{Gr} T. Damm, L. Gr\"une,  M. Stieler,  K. Worthmann, An exponential turnpike theorem 
for dissipative discrete time optimal control problems. SIAM J. Control Optim. 52  (2014), 1935-1957.

\bibitem{Gugat1} M. Gugat, Optimal energy control in finite time by varying the length of the string. SIAM J. Control
Optim. 46 (2007), 1705-1725.

\bibitem{GTZ} M. Gugat, E. Tr\'elat, E. Zuazua, Optimal Neumann control for the 1D wave equation: Finite horizon, infinite horizon, boundary tracking terms and the turnpike property. Systems and Control Letters 90 (2016), 61-70.



\bibitem{HS}
A. Henrot,  J. Sokolowski, 
A shape optimization problem for the heat equation. Optimal Control (Gainesville, FL, 1997), 204-223, Appl. Optim., 15, Kluwer Acad. Publ., Dordrecht, 1998.


\bibitem{HHS}
M. Hayouni, A. Henrot, N. Samouh,  On the Bernoulli free boundary problem
and related shape optimization problems. Interfaces \& Free Bound. 3 (2001), 1-13.

\bibitem{HP1}
A. Henrot, M. Pierre, Variation et Optimisation de Formes. une analyse g\'eom\'etrique., Math\'ematiques et Applications. 48 (2005).

\bibitem{LR}
W. B. Liu, J. E. Rubio,  Optimal shape design for systems governed by variational inequalities. Part 1: Existence theory for the elliptic case, Part 2: Existence theory for the evolution case. J. Optim. Theory  Appl. 69  (1991), 351-371, 373-396.


\bibitem{Ivan}
I. Moyano, Controllability of a 2D quantum particle in a time-varying disc with radial data.
https://hal.archives-ouvertes.fr/hal-01405624v1.


\bibitem{P1}
O. Pironneau,  Optimal Shape Design for Elliptic Systems. Springer-Verlag,
1984.

 \bibitem{PZ1} A. Porretta, E. Zuazua, Long time versus steady state optimal control. 
SIAM J. Control and Optim. 51 (2013), 4242-4273.

\bibitem{PZ2} A. Porretta, E. Zuazua, Remarks on long time versus steady state optimal control. Springer-INdAM, Mathematical Paradigms of Climate Science, P. M. Cannarsa et al. eds, Springer International Publishing Switzerland, 2016, 67-89.

\bibitem{PTZ2}
Y. Privat, E. Tr\'elat, E. Zuazua,
Optimal shape and location of sensors for parabolic equations with random initial data. Arch. Ration. Mech. Anal. 216 (2015), 921-981.

\bibitem{PTZ1}Y. Privat, E. Tr\'elat, E. Zuazua, Optimal observability of the multi-dimensional wave and Schr\"odinger equations in quantum ergodic domains. J. Eur. Math. Soc. 18 (2016), 1043-1111.

\bibitem{PTZ3}
Y. Privat, E. Tr\'elat, E. Zuazua,
Randomised observation, control and stabilization of waves [Based on the plenary lecture presented at the 86th Annual GAMM Conference, Lecce, Italy, March 24, 2015]. ZAMM Z. Angew. Math. Mech. 96 (2016),  538-549. 


\bibitem{S11}
J. Simon, Differentiation with respect to the domain in boundary value problems. Numer. Funct. and Optimiz.,
 2 (1980), 649-687.

\bibitem{asz}
J. Sokolowski, J. P.  Zolesio, Introduction to Shape Optimization. Springer Berlin Heidelberg, 
1992.

\bibitem{SZ}
J. Sokolowski, A. Zochowski, On the topological derivative in shape optimization. SIAM J. Control
Optim. 37 (1999), 1251-1272.

\bibitem{S1} V. \u{S}ver\'{a}k, On optimal shape design. J. Math. Pures Appl. 72 (1993), 537-551.

\bibitem{Tartar}
L. Tartar, An introduction to the homogenization method in optimal design. Optimal shape design. Springer Berlin Heidelberg, 2000, 47-156.

\bibitem{TZ1} E. Tr\'elat, E. Zuazua, The turnpike property in finite-dimensional nonlinear optimal control.
 J. Differential Equations 258 (2015), 81-114.

\bibitem{TZZ1} E. Tr\'elat, C. Zhang, E. Zuazua, Integral and measure-turnpike property for infinite-dimensional optimal control problems. Preprint.

\bibitem{TZZ2} E. Tr\'elat, C. Zhang, E. Zuazua,  Steady-state and periodic exponential turnpike property for optimal control problems in Hilbert spaces. https://arxiv.org/abs/1610.01912.

\bibitem{Marius} M. Tucsnak, G. Weiss, Observation and Control for Operator
Semigroups, Springer, 2009.


\bibitem{Z1}
A. J. Zaslavski, Turnpike Properties in the Calculus of Variations and Optimal Control. Vol. 80. Springer, 2006.

\bibitem{Z2}
A. J. Zaslavski, Turnpike Theory of Continuous-time Linear Optimal Control Problems. Springer Optimization and Its Applications, 104. Springer, Cham, 2015.

\end{thebibliography}
\end{document}